\documentclass[11pt]{article}
   
\usepackage{lineno}

\usepackage{caption}
\usepackage{color}
\usepackage{longtable}
\usepackage{subcaption}
\usepackage{afterpage}
\usepackage{amsmath}
\usepackage{amssymb}
\usepackage{graphicx}
\usepackage{parskip}
\usepackage{float}
\usepackage{multirow}
\usepackage{booktabs}
\usepackage{pgfplots}
\usepackage{tabularx}
\usepackage{hyperref}
\usepackage{algorithm}
\usepackage{algpseudocode}
\usepackage[doublespacing]{setspace}

\usepackage{layout}
\usepackage{makecell}
\usepackage{mwe}

\usepgfplotslibrary{statistics}
\usepackage{booktabs}

\usepackage{pgfplots}
\pgfplotsset{compat=1.11,
        /pgfplots/ybar legend/.style={
        /pgfplots/legend image code/.code={%
        \draw[##1,/tikz/.cd,bar width=3pt,yshift=-0.2em,bar shift=0pt]
                plot coordinates {(0cm,0.8em)};},
},}
\usepackage{tikz}
\usetikzlibrary{shapes.geometric, arrows}
\tikzstyle{startstop} = [rectangle, rounded corners, minimum width=2cm, minimum height=1cm,text centered, draw=black, fill=red!30, text width=2cm]  
\tikzstyle{io} = [trapezium, trapezium left angle=70, trapezium right angle=110, minimum width=2cm, minimum height=1cm, text centered, draw=black, fill=blue!30, text width=2cm]
\tikzstyle{process} = [rectangle, minimum width=2cm, minimum height=1cm, text centered, draw=black, fill=orange!30, text width=2cm]
\tikzstyle{decision} = [diamond, minimum width=2cm, minimum height=1cm, text centered, draw=black, fill=green!30, text width=2cm]
\tikzstyle{arrow} = [thick,->,>=stealth]
\usetikzlibrary{calc}

\usepackage[thmmarks]{ntheorem}
\usepackage{amssymb}

\theoremstyle{plain}

\theoremstyle{nonumberplain}
\theoremsymbol{\ensuremath{\Box}}

\makeatletter
\newtheoremstyle{nonumberplainnobrackets}%
  {\item[\theorem@headerfont\hskip\labelsep ##1\theorem@separator]}%
  {\item[\theorem@headerfont\hskip \labelsep ##1\ ##3\theorem@separator]}
\makeatother

\providecommand{\keywords}[1]
{	
  \textbf{\textit{Keywords---}} #1
}
\theoremstyle{nonumberplainnobrackets}
\theoremsymbol{\ensuremath{\blacksquare}}

\usepackage{algorithm}

\DeclareMathAlphabet{\mathcal}{OMS}{cmsy}{m}{n}
\SetMathAlphabet{\mathcal}{bold}{OMS}{cmsy}{b}{n}

\usepackage{url}
\usepackage[round,longnamesfirst]{natbib}
\usepackage[final]{changes}

\title{A mixed-integer programming model for identifying intuitive ambulance dispatching policies}
 \author{\textbf{Laura A.~Albert} \\
 	{\small \vspace{-6pt}Department of Industrial and Systems Engineering} \\
 	{\small \vspace{-6pt}University of Wisconsin-Madison} \\
 	{\small \vspace{-6pt}1513 University Avenue, Madison, WI 53706} \\
 	{\small \vspace{-6pt}\url{laura@engr.wisc.edu},  Twitter: \href{http://twitter.com/lauraalbertphd}{@lauraalbertphd} }  \\
 	{\small \vspace{-6pt}ORCID: 0000-0001-7079-4473 } }
\date{October 2022}

\begin{document}

\maketitle

\normalsize

\abstract{Markov decision process models and algorithms can be used to identify optimal policies for dispatching ambulances to spatially distributed customers, where the optimal policies indicate the ambulance to dispatch to each customer type in each state. Since the optimal solutions are dependent on Markov state variables, they may not always correspond to a simple set of rules when implementing the policies in practice.  Restricted policies that conform to a priority list for each type of customer may be desirable for use in practice, since such policies are transparent, explainable, and easy to implement. A priority list policy is an ordered list of ambulances that indicates the preferred order to dispatch the ambulances to a customer type subject to ambulance availability.  This paper proposes a constrained Markov decision process model for identifying optimal priority list policies that is formulated as a mixed integer programming model, does not extend the Markov state space, and can be solved using standard algorithms.  A series of computational examples illustrate the \added{benefit of intuitive policies}. The optimal mixed integer programming solutions to the computational examples have objective function values that are close to those of the unrestricted model and are superior to those of heuristics. }


\section{Introduction} 
This paper is concerned with finding intuitive dispatching policies for sending ambulances to spatially-located customers that are straightforward to implement. Optimal solutions to the dispatching problem can be
identified using Markov decision process (MDP) models and
algorithms, where a solution to the MDP indicates the optimal ambulance to send in each Markov state. Therefore, the optimal policy is a collection of a
large number of decisions that may be difficult to succinctly
summarize and implement in practice. As a result, full state-dependent model solutions may send different ambulances to the same type of patient in two seemingly equivalent situations corresponding to different states. In this paper, we seek policies that are optimal in some sense but are not fully state dependent. 
These optimal \emph{restricted policies} are more explainable than the optimal unrestricted policies that are fully state-dependent. We formulate new mixed integer programming modela to identify restricted policies in an MDP that conform to a priority list, and we compare the resulting policies to fully state-dependent MDP policies. To our knowledge, there are no papers that have studied this issue.

A \emph{priority list} policy is a policy that, for each customer
type, dispatches ambulances in a given order based on ambulance
availability. Priority list policies are advantageous in that they can be implemented  by computer software and explained by human decision-makers using a simple set of rules or a decision tree. An example of a priority list policy is to dispatch
ambulance 1 first (i.e., in all MDP states where ambulance 1 is an available
choice), ambulance 2 only if ambulance 1 is busy, ambulance 3 if ambulances 1
and 2 are busy, and so on. Each of these choices reflects many states in the MDP. In contrast, the optimal unrestricted MDP policy may not consistently send the same ambulance. For example, instead of sending ambulance 1 first if it is available, the first ambulance to send may depend on the available ambulances, e.g,. it may send ambulance 1 if all ambulances are available, ambulance 2 if ambulance 3 is busy (and ambulance 1 is free), ambulance 3 if ambulance 4 is busy, and so on. The set of priority list policies is a
subset of the set of all possible MDP policies, yet limiting the type of
policies considered within the MDP model is not straightforward,
since the MDP state space is not structured to identify restricted
policies that require decisions to be the same in distinct states
unless the state space is expanded to account for this restriction. 

The use of intuitive policies in emergency dispatching in practice has important implications for situations when life-saving decisions must be made in an extremely short amount of time. Policies that are simple to implement are advantageous when dispatch decisions are made under duress, such as during the theater of military operations or during emergency response operations during catastrophic events. In these applications, the full set of policies may not be feasible or desirable to implement.
Explainability is widely acknowledged feature for the practical deployment of artificial intelligence tools and decision support systems \citep{arrieta2020explainable}. Explaining algorithmic decisions has also become a legal requirement in the European Union \citep{goodman2017european}.
Explainable policies are attractive to public safety leaders who often must report performance metrics to the public and practical for dispatchers who have the discretion regarding which vehicles to send to calls for service. For example, New York City Police Department’s Domain Awareness System uses logic-based and threshold-based alerts developed based on the feedback of police officers to ensure their explainability and to facilitate the usage of the decision support system in the field \citep{levine2017new}. 

Several recent papers that study optimal dispatching using MDPs \citep{jagtenberg2017dynamic} and approximate dynamic programming \citep{schmid2012solving,robbins2020approximate,jenkins2021approximate} do not require that dispatching policies conform to priority lists.
However, priority list policies are widely used in more routine operations, such as those faced by emergency medical services (EMS) \citep{Ansari-etal-2012}. 
Simple heuristics used in practice, such as the heuristic of sending the closest
available ambulance relative to a customer, follow a priority list.
Models that compute exact or approximate dispatch probabilities
in spatial queueing systems often implicitly assume
that the queueing preferences correspond to a priority list
\citep{Larson74,Jarvis85,Budge-etal-09,rautenstrauss2021ambulance}.
Other models for locating ambulances
\citep{Erkut-etal-08,SilvaSerra2008,toro2015reducing,belanger2020recursive} or jointly locating
and dispatching ambulances
\citep*{Ansari-etal-2012,ToroDiaz-etal-2012} often assume a
priority list policy for dispatch. The existing papers that identify priority list policies do not compare the priority list policies to the best overall policies that may not conform to a priority list. This paper addresses this gap \added{by proposing a methodology to identify optimal, intuitive policies}.

Relocating ambulances upon completion of service relates to priority lists. Many emergency medical services use \emph{compliance tables} to locate and relocate ambulances when they finish servicing a call \citep*{Alanis-etal-2012}. In a compliance table, the location and relocation decisions are based solely on the number of available ambulances, which makes such policies simple to implement. The sets of ambulance locations are nested, which is an important feature of compliance tables \citep*{van2017compliance}. 
A compliance table can be created through integer programming \citep*{gendreau2006maximal}, Hypercube queueing models \citep*{van2017compliance}, and Markov models \citep*{Alanis-etal-2012}, and they yield the desired set of ambulance locations to use given the number of available ambulances. If the ambulances are at these locations, then the system is said to be compliant. 
Similarly, \cite{Dimitrov-etal-2011} examine how to locate sensors in a network where the number of sensors is uncertain such that the sensor locations are nested like compliance tables. This requirement makes sensor deployment more practical. The priority list requirement does
not substantially reduce the objective function values in the case of sensor placement \citep*{Dimitrov-etal-2011}, a finding similar to what we observe in the computational examples in Section 4.  

Priority list policies can be identified by constrained MDP models \citep{Altman99}.
Enforcing a priority list policy can also be viewed as a type of MDP design
problem, a one-time set of decisions that is made prior to the control process. Then, the MDP optimally controls the process given the design choices. MDP design problems can be modeled as linear or mixed integer programs that do not expand the state space, which is similar to the approach we take in this paper.  \cite{mclay2013equity} formulate a constrained MDP using linear programming to identify equitable ambulance dispatching policies.  \cite{Dimitrov-Morton2009} consider costs associated with the design variables and simple relationships between the design variables and linear MDP variables. Neither of these papers imposes a particular structure on the policies as we do in this paper. Instead, they illustrate how to model resource allocation or equity side constraints.
In contrast, this paper considers priority list decision decisions that have no associated costs but rather are used to restrict the structure of the optimal policies.


This paper proposes a constrained MDP model formulated as a mixed integer programming (MIP) model to
identify optimal priority list policies for dispatching ambulances to spatially located customers. The approach extends the standard linear programming formulation for
unconstrained MDP models where no structure is imposed on the
optimal policy. In the proposed model, priority list policies are
enforced by introducing binary variables to capture the priority lists as well as an additional set of constraints. Each
binary variable corresponds to an ambulance and its priority to send to
a customer in the priority list. 
New constraints link the binary priority list variables to the linear dispatching variables to shape
the structure of the optimal policies, requiring some actions to be
the same in different Markov states, and the model solutions represent restricted optimal policies that are not
fully state-dependent. The resulting formulation does not expand the Markov
state space, and the MIP model formulation can be solved using standard
mixed integer programming algorithms. The proposed model is developed for a model for dispatching
ambulances to prioritized customers in an average cost, undiscounted,
infinite horizon MDP
\citep{McLay-Mayorga2012dispatch}.
The MIP model focuses on the decision context
of dispatching ambulances (ambulances) to prioritized customers
(patients). The base model only allows non-idling policies and has a zero-length queue, leading to customers that are ``lost''
only when all ambulances are busy. The base model is then extended to allow for idling policies to some customer types.  A computational example is
analyzed to shed light on when the unrestricted optimal MDP policies
are (and are not) priority list policies, how the unrestricted and
restricted policies differ, and how the priority list policies
improve upon heuristic policies that also conform to priority lists. \added{While the ultimate goal is to identify intuitive policies and to demonstrate their benefit, a discussion of how the approach scales is also discussed.}

This paper is organized as follows.  The linear programming
model for the unrestricted model is presented in Section 2. The
proposed constrained MDP models \added{for identifying optimal, restricted policies} are formulated in Section 3
using mixed integer programming models. A computational example
is presented and analyzed in Section 4. Concluding remarks are
presented in Section 5.

\section{\added{Unrestricted model formulation}}
This section presents an MDP model for
assigning distinguishable ambulances to different types of customers.
No structure is imposed on the MDP solutions, which thus yield
the unrestricted solutions. The objective is to optimally determine
which ambulance to dispatch to arriving customers to maximize
the average utility per stage in an average cost, undiscounted, infinite horizon MDP. Following the introduction of the MDP
model, the unrestricted MDP is formulated as a linear programming
model. In the next section, we extend the linear programming model
to require that the solutions follow a priority list structure.

In the MDP models, there is a set of customer types $W^+$ who arrive according to a Poisson process with rate $\lambda_w, w \in W^+$. As soon as a customer arrives to the system at one of the $n$ locations, its type becomes known, and one of the $m$ available ambulances is dispatched to it if one is available. This implies the use of non-idling policies, which is reasonable since the model is motivated by EMS systems where providing service to all customers is a major goal. However, some research has considered idling policies \citep{yoon2018expected,dubois2021dispatching}. To simplify the model presentation, we also include a ``null'' customer to represent MDP transitions when a customer does not arrive and an ambulance is not assigned. Let $W$ capture the entire set of customer types, including the null customer. Table \ref{table:notation} summarizes the MDP notation.

Any ambulance could be assigned to any customer if the ambulance is available. Ambulances start service at a home station and return to the same station after servicing a customer, and multiple ambulances could share the same home station. This assumption reflects a static deployment of ambulances that are not dynamically relocated after servicing customers. 
Service is non-preemptive, which is a common restriction in practice. Service times depend on the customer location and are exponentially distributed, which is not realistic yet is required to keep the state space manageable. Finally, we assume a zero-length queue, which implies that customers are ``lost'' when they arrive when all ambulances are busy. In emergency medical service systems, for example, the lost
customers are often treated by service providers other than
paramedics and emergency medical technicians, such as fire engines. Research suggests that the exponentially distributed service time and zero-length queue assumptions are reasonable when the analytical results are evaluated using simulation \citep{bandara2014priority,McLay-Mayorga2012dispatch}.
In Sections 3.1, we lift the assumption that an ambulance must always be sent if one is available.

\begin{table}
\singlespacing
\caption{Markov decision process model notation}\label{table:notation}
\begin{tabular}{l c l}
\hline 
$m$ &=&  the number of ambulances, each at a fixed base while not in service \\
$n$ &=&  set of customer locations \\
$S$ &=&  set of states, where each state corresponds to the
  configuration of ambulance \\
$W$ &=&  the set of customer types, including a null   customer. Let $W^+ \subseteq W$ capture \\
&& the set of non-null customer arrivals \\
$\lambda_w$ &=&  the arrival rate for customers of type $w \in W^+$ with overall arrival rate \\
&& $\lambda = \sum_{w \in W^+}\lambda_w$ \\
$\mu_{ij}$ &=&  the service rate (distributed exponentially) when ambulance $j$ is dispatched \\
&& to a customer of any type located at $i=1,...,n$ \\
$X(s,w)$ &=&  the set of actions available in state
  $s \in S$ when customer type $w \in W$ arrives,\\
  && including a null action. Let $X^+(s,w) \subseteq X(s,w)$ capture the set of \\
  && ambulances available to send to customer types $w \in W^+$. \\
$Pr(s',w' | s,w,j)$ &=&  the conditional probability that the
  system moves from state $s \in S$ with \\
  && customer type $w \in W$ to
  state $s' \in S$ with customer type $w' \in W$ \\
  && when action $j \in X(s,w)$ is selected \\
$r(s,w,j)$ &=& the reward when ambulance $j$ is dispatched to customer type $w$ in
state $s$ \\
\hline
\\
\end{tabular}
\end{table}

The state of the system is defined as a vector of ambulance locations,
$S=(s_1, s_2, ..., s_m)$. If ambulance $j$ is in location 0 (i.e., $s_j=0$), it is
available; otherwise it is busy servicing a customer at location $i =1,2,...,n$ (i.e., $s_j=i$).
The customer types $W^+$ reflect combinations of location and triage level. The computational example in Section 4 admits high-level and low-level customers $h \in \{H,L\}$ that arrive at one of $n$ locations. This leads to $|W| = 2n$ customer types not including null customer types and $|S|=(n+1)^m$. Note that if service times are identical across locations, as often assumed in the literature (see \cite*{Budge-etal-09}), then the state space can be substantially reduced to $|S|=2^m$. The larger state space is used here to highlight the value of \added{intuitive policies, which we explore in greater detail in the case study in in Section 4}.

Let $X(s,w)$ denote the set of available actions in state $s$ given
that a customer with type $w$ arrives. The set of actions captures the available ambulances with $X(s,w) \subseteq \{1,2,...,m\}$,  or if no ambulances are available, a null ambulance  (ambulance 0).  When the triggering event is not a customer arrival, the null actions are
 either to move an ambulance completing service to its home location or to do nothing.
 Let $X^+(s,w) \subseteq X(s,w)$ capture the set of non-null actions associated with sending an ambulance to a customer. The rewards $r(s,w,j)$ reflect the
utility of the action selected, which depends on the customer type
as well as the ambulance that responds. Typically, the utility in emergency medical dispatch reflects the proportion of calls responded to within a fixed timeframe (e.g., nine minutes), which reflect the travel times between ambulance $j$ and location $i$  \citep*{Erkut-etal-08,McLay-Mayorga2012dispatch}, although the utilities can be defined to capture other criteria as well \citep{yoon2018expected,yoon2021dynamic}. A reward of zero is associated
with assigning a null ambulance or a ``lost'' customer. The transition
probabilities capture (1) the busy ambulances completing service and
moved to location 0, and (2) a customer arriving, which requires
that an ambulance is dispatched to the customer if one is available.
Under these modeling assumptions, the MDP is recurrent, since each
deterministic stationary policy leads to a single class where any
ambulance could be assigned to any customer type, even if it is the
last possible choice.

The model is formulated as an equivalent discrete time MDP
using uniformization. 
The maximum rate of transitions is
determined as $\gamma = \lambda + \sum_{j=1}^m \beta_j$ where
$\beta_j = \max_{i=1,2,...,n} \{ (1 / \mu_{ij}) \}$. The
transition probabilities $Pr(s',w' | s,w,j)$ are computed as
follows. The probability that ambulance $j$ becomes available if
it is serving a customer of any type at location $i$ is $(\gamma
\mu_{ij})^{-1}$. This causes a transition from state $s$ with
$s_j=i$ to the corresponding state $s'$ with $s'_j=0$.
Assigning ambulance $j$ to an arriving customer type $w$ at location $i$ in state $s$ with $s_j=0$ causes a transition to the corresponding state $s'$ with $s'_j=i$ with probability $\lambda_w/ \gamma$. Finally, the probability of staying in state $s$ is $\left( 1 - \lambda/\gamma - \sum_{j=1}^m
\sum_{i>0: s_j = i} \frac{1}{\gamma \mu_{ij}}\right)$, or
$\left( 1 - \sum_{j=1}^m \sum_{i>0: s_j = i} \frac{1}{\gamma
\mu_{ij}}\right)$ if no ambulances are available in $s$. Let $Pr(s',w' | s,w,j)$ capture the resulting transition probabilities that the system moves from state $s \in S$ with customer type $w \in W$ to state $s' \in S$ with customer type $w' \in W$ when action $j \in X(s,w)$ is selected, which implicitly captures customer arrivals and service completions. The
reader is referred to \cite{McLay-Mayorga2012dispatch} for
additional details.

The optimal policy for the unconstrained MDP can be identified by
solving the linear programming (LP) formulation and yields a fully state-dependent dispatching policy. As noted earlier, the system is described by the state (the configuration of ambulances) as well as the type of customer that has arrived ($w$). An action must be
selected when an event occurs, where an event could be a new
customer, a service completion, or a null event due to
uniformization.  The variables in the LP model $y(s,
w, j)$, $s \in S, w \in W, j=1,2,...,m$, represent the
proportion of time the system is in state $s$ when customer $w$ has
arrived and action $j$ is selected.

First, the unrestricted MDP model (U) is formulated as a LP
model.
\begin{align}
\label{linprogmodel1}\mathrm{Model\ U:} \ & Z^U = \max \sum_{s \in S, w \in W, j \in X(s,w)}
\!\!\!\!\!\!\!\!\!\!\!\!r(s,w,j)\, y(s,w,j) & \\
\label{linprogmodel2} \mathrm{subject\ to}&
\sum_{j' \in X(s',w')} \!\!\!\!\!\!
y(s',w',j') - \!\!\!\!\!\!\!\!\!\sum_{s \in S, w \in W, j \in
X(s,w)} \!\!\!\!\!\!\!\!\!\!\!\!\!\!\!\! Pr(s', w' |s,w,j)\, y(s,w,j)=
0, s'
\in S, w' \in W & \\
\label{linprogmodel3} & \sum_{s \in S, w \in W, j \in X(s,w)}
\!\!\!\!\!\!\!\!\!\!\!\! y(s,w,j)=1 & \\
\label{linprogmodel4} & y(s,w,j) \ge 0, s \in S, w \in W, j \in X(s,w). &
\end{align}
In Model U, the objective function \eqref{linprogmodel1} reflects the average reward per
stage. The first set of constraints \eqref{linprogmodel2} balances the flow in and out of
each state and customer type. The second set of constraints \eqref{linprogmodel3} requires
that all of the actions sum to one. The final set of constraints \eqref{linprogmodel4}
requires that the variables are non-negative. Note that since there
is a finite number of states and bounded rewards, the optimal policy
is Markovian and deterministic \citep{Puterman94}. 

\section{\added{Intuitive model formulation}}\label{section:mip} 
In this section we present a mixed
integer programming model formulation to identify optimal priority
list policies assocated with the MDP introduced in the previous section. To do so, we must identify a rank ordering of the $m$ ambulances for each customer type $w \in W^+$. Note that if two ambulances are located at the same station, both ambulances are ordered in the priority list.
There are two sets of variables. The first set of variables $y(s,w,j)$, $s \in S, w \in W,
j=1,2,...,m$, is defined in the same way as in Model U. Recall that
we are concerned with the structure of how (non-null) ambulances are
dispatched to (non-null) customers; other types of decisions selected
in the MDP are not a principle concern here. The second set of
variables comprise the binary priority list design variables $z$,
where $z_{w,j}^p$ is 1 if ambulance $j$ is the $p$th preferred ambulance
to dispatch when customer $w \in W^+$ arrives and 0 otherwise, with
$|z| = |W^+|m^2$. Define $P_{s,w} = m-|X(s,w)|+1$, which is 1 when
all of the (non-null) ambulances can be dispatched to customer type $w$
and $m$ when one non-null ambulance can be dispatched to customer type
$w \in W^+$ in state  $s \in S$. Here $P_{s,w}$ indicates the
largest priority of an ambulance in the list that could be dispatched in state $s
\in S$ to customer $w \in W^+$.

The mixed integer programming model (Model PL) for priority list policies is formulated by
adding constraints (\ref{mipmodel1}) -- (\ref{mipmodel4}) to
Model U (given by (\ref{linprogmodel1}) --
(\ref{linprogmodel4})) with objective function value $Z^{PL}$.
\begin{align}
& \label{mipmodel1}  \sum_{p=1}^m z_{w,j}^p = 1, w \in W^+,
j=1,2,...,m & \\
& \label{mipmodel2}  \sum_{j=1}^m z_{w,j}^p = 1, w \in W^+, p=1,2,...,m & \\
& \label{mipmodel3} y(s,w,j) \le 1 - \!\!\!\!\! \sum_{j' \in X^+(s,w) \setminus j' }
\!\!\!\!\! z_{w,j'}^p + \sum_{p'=1}^{p-1} z_{w,j}^{p'},  \ s \in S,w \in W^+, j \in X^+(s,w),
p=1,2,...,P_{s,w} & \\
& \label{mipmodel4}  z_{w,j}^p \in \{0,1\}, w \in W^+, \ j=1,2,...,m, p=1,2,...,m.&  
\end{align}
The first and second sets of constraints, (\ref{mipmodel1}) and
(\ref{mipmodel2}), require every (non-null) ambulance to be
represented by exactly one priority and every priority to be
assigned exactly one ambulance.
The third set of constraints (\ref{mipmodel3}) links the dispatch variables $y(s,w,j)$ to the priority lists to ensure the MDP decisions are harmonized across states. It does so by using the priority list dispatch variables $z_{w,j}^p$ to set variable upper bounds for each of the linear dispatching variables $y(s,w,j)$. This set of constraints can reduce the variable upper bound for a $y(s,w,j)$ variable from one to zero---thereby removing the action---if another action is available with a higher priority based on the values of the priority list dispatch variables $z_{w,j}^p$. The
variable upper bound is one if the action under consideration is selected prior
to other allowable actions that are a lower priority (shown by the
last term in (\ref{mipmodel3})). This set of constraints limits the
allowable actions to one action per state based on the priority list
variables. Note that several such upper bounds may be set for each linear
variable $y(s,w,j)$, and hence, an action is removed if just one of
these upper bounds is set to zero. The last set of constraints (\ref{mipmodel4}) ensures that the priority list design variables are binary. 

Constraints (\ref{mipmodel3}) are illustrated through the
following example with $m=3$ ambulances. Assume that the
priority list for customer $w$ is to first send ambulance 3, then
ambulance 2, and ambulance 1 last. Then,
$z_{w,3}^1=z_{w,2}^2 = z_{w,1}^3=1$, and $z_{w,j}^p=0$ for all other
combinations of $j$ and $p$ for customer type $w$. Consider state $s'$, in
which all three ambulances are available, with $X(s',w)=\{1,2,3\}$ and
$P_{s',w}=1$ (the most preferred ambulance will be sent, which is ambulance 3). Then
(\ref{mipmodel3}) results in the following three constraints being
added to Model PL:
\begin{eqnarray}
\nonumber p=1,& j=1: & y(s',w,1) \le 1 - z_{w,2}^1 - z_{w,3}^1 = 0\\
\nonumber p=1,& j=2: & y(s',w,2) \le 1 - z_{w,1}^1 - z_{w,3}^1 = 0\\
\nonumber p=1,& j=3: & y(s',w,3) \le 1 - z_{w,1}^1 - z_{w,2}^1 = 1
\end{eqnarray}
These three constraints provide upper bounds for the three variables corresponding to the three available actions. The upper bounds for ambulances 1 and 2 are set to zero,
these actions are not allowed in this state. Therefore, only action
$j=3$ is allowed, which is the first priority action.

Consider the same customer type $w$ in a different state $s''$
with ambulances 1 and 2 available. Then, $X(s'',w)=\{1,2\}$ and
$P_{s'',w}=2$. Then, $z_{w,3}^1=z_{w,2}^2=z_{w,1}^3=1$, and
$z_{w,j}^p=0$ otherwise (the same as before), and
(\ref{mipmodel3}) results in the following four constraints
being added to Model PL:
\begin{eqnarray}
\nonumber p=1,& j=1: & y(s'',w,1) \le 1 - z_{w,2}^1  = 1\\
\nonumber p=2,& j=1: & y(s'',w,1) \le 1 - z_{w,2}^{2} + z_{w,1}^1 = 0\\
\nonumber p=1,& j=2: & y(s'',w,2) \le 1 - z_{w,1}^1  = 1\\
\nonumber p=2,& j=2: & y(s'',w,2) \le 1 - z_{w,1}^{2} + z_{w,2}^1 = 1
\end{eqnarray}
The first two constraints correspond to ambulance 1, and the last two
constraints correspond to ambulance 2. The second constraint disallows
ambulance 1 from being dispatched by setting its upper bound to zero.
Both upper bounds for $y(s'',w,2)$ are 1, which means that these
bounds do not limit how $y(s'',w,2)$ can be set, thus leading to
this action being selected.


When a single action is available in a state, constraint set (\ref{mipmodel3})
degenerates to trivial upper bound of 1 for the linear variable
corresponding to the one available action.
Therefore, when there is only one action available in for a
state-customer pair, then variable upper bound constraints do
not need to be added.

The MIP formulations suggest that the priority list requirement can
be viewed as a type of design decision, where a one time priority
list is first designed that determines the states during the MDP
process. 
The MIP models could simultaneously consider other design decisions. This could be achieved by adding a set of budget constraints to the model formulation that enforce multiple knapsack constraints on the set of priority list design variables. In this case, knapsack $k$ (of $K$ total knapsacks) has capacity $b_k$ and item $z_{w,j}^p$ has weights $a_{w,j}^{p,k}$, $k=1,2,...,K$, yielding the set of constraints 
\begin{displaymath}
\sum_{w \in W_0}\sum_{j=1}^m\sum_{p=1}^m a_{w,j}^{p,k} z_{w,j}^p \le
b_k,\ k=1,2,...,K.
\end{displaymath}
However, these constraints may not be appropriate for the emergency
dispatch applications considered in this paper.

\subsection{A MIP model that allows idling}\label{section:mipidling}
In this section, we formulate Model PLI that lifts the non-idling policy assumption for some customer types to allow for the action space to allow for some customers to be ``lost'' even when there are available ambulances. It is assumed that customers that are lost are served by a neighboring region through mutual aid. We partition the set of (non-null) customer types into low-level ($L$)
and high-level ($H$) subsets, i.e., $W^+ = W^+_L \cup W^+_H$, with $W^+_L \cap W^+_H = \emptyset$. 
To maintain a recurrent MDP, we assume that an idling policy can only be used only for the low-level customers.
The additional idling action for low-level customers in $W^+_L$ is captured by ambulance 0. Therefore, there are a total of $m+1$ ``ambulances'' that are ordered for each low-level customer type.  This new ``ambulance'' corresponding to the idling action does not change the state space, and it is modeled as an ambulance that is always available for low-level customers in the action sets $X(s,w)$. The idling action is available in all MDP states, and therefore, the actions that appear in the priority list after the idling ambulance 0 are never selected by the MDP, although they appear in the priority lists.  Based
on these assumptions, define $P^L_{s,w}=(m+1)-|X(s,w)|+1$ and
$P^H_{s,w}=(m)-|X(s,w)|+1$. Model PLI is formally stated as (\ref{pliq1}) -- (\ref{pliq2}).
\begin{eqnarray}
\label{pliq1}\mathrm{Model\ PLI} &Z^{PLI} = \max & \sum_{s \in S, w \in W, j
\in
X(s,w)} \!\!\!\!\!\!\!\!\!\!\!\!r(s,w,j)\, y(s,w,j) \\
\nonumber \mathrm{subject\ to}&& \!\!\!\!\!\!\!\!\!\!\!\!\sum_{j' \in
X(s',w')} \!\!\!\!\!\! y(s',w',j') - \!\!\!\!\!\!\!\!\!\!\!\!\!\!\!\!\!\!\sum_{s \in S, w \in W, j \in X(s,w)} \!\!\!\!\!\!\!\!\!\!\!\!\!\!\!\! Pr(s', w'|s,w,j)\,
y(s,w,j)= 0, \\
&& \hspace{2.0in} \mathrm{for}\ s' \in S, w' \in W\\
 && \sum_{s \in S, w \in W, j \in X(s,w)}
\!\!\!\!\!\!\!\!\!\!\!\! y(s,w,j)=1\\
 && y(s,w,j) \ge 0, s \in S, w \in W,j \in X(s,w).
\end{eqnarray}
The objectives and constraints thus far are analogous to \eqref{linprogmodel1} -- \eqref{linprogmodel4}. The next five sets of constraints set the priority list variables
for the low-level customers that allow for idling.
\begin{eqnarray}
 && \sum_{p=1}^{m+1} z_{w,j}^p = 1, w \in W_L^+, j=0,1,...,m\\
 && \sum_{j=0}^{m} z_{w,j}^p = 1, w \in W_L^+, p=1,2,...,m+1\\
\nonumber && y(s,w,j) \le 1 - \sum_{j' \in X^+(s,w) \setminus j'}
 z_{w,j'}^p + \sum_{p'=1}^{p-1} z_{w,j}^{p'},\\
 && \hspace{1.0in} \mathrm{for}\ s \in S, w \in W_L^+, j \in
X^+(s,w), p=1,2,...,P^L_{s,w} \\
  && z_{w,j}^p \in \{0,1\}, w \in W_L^+, j=0,1,...,m, p=1,2,...,m+1.
\end{eqnarray}
The final five sets of constraints set the priority list
variables for the high-level customers that do not allow idling. These are the same as the analogous constraints in Model PL.
\begin{eqnarray}
 && \sum_{p=1}^{m} z_{w,j}^p = 1, w \in W_H^+, j=1,...,m\\
 && \sum_{j=1}^{m} z_{w,j}^p = 1, w \in W_H^+, p=1,2,...,m\\
\nonumber && y(s,w,j) \le 1 - \sum_{j' \in X^+(s,w) \setminus j'}
 z_{w,j'}^p + \sum_{p'=1}^{p-1} z_{w,j}^{p'},\\
 && \hspace{1.0in} \mathrm{for}\ s \in S, w \in W_H^+, j \in
X^+(s,w), p=1,2,...,P^H_{s,w} \\
\label{pliq2} && z_{w,j}^p \in \{0,1\}, w \in W_H^+, j=1,...,m, p=1,2,...,m.
\end{eqnarray}

It is worth noting that in many public safety systems, such as
emergency medical services, idling policies are not typically used,
since all customers are expected to be served immediately. However, some recent research has considered the use of idling policies for low-level patients \citep{yoon2018expected}. Additionally, low-level patients are sometimes not immediately served during natural disasters or mass casualty events when the number of patients who require service overwhelm the system \citep{dubois2021dispatching}.

\section{Computational examples}
In this section, we illustrate the solutions to Models PL and PLI,
and we compare the resulting priority list policies from Model U and from a heuristic. The goal is to shed light on how
restricted policy changes the underlying decisions and to
understand when the priority list solutions improve upon those
obtained by using a heuristic. The heuristic that we use is the
closest ambulance policy, which dispatches the closest available ambulance relative to a customer and conforms to a priority
list. \added{Then, we discuss how the approach scales.}

The examples consider the following geographic regions with $m=4$ ambulances. Each of the five region maps shown below is composed of $n=4$ locations where customers can arrive, and an ambulance is located in the center of each customer location, both of which are numbered from 1 to 4.

\setlength{\unitlength}{0.5cm}
\begin{picture}(4,1.5)
\put(0,0.3){Region R1} \put(4,0){\line(0,1){1}}
\put(4,0){\line(1,0){4}} \put(4,1){\line(1,0){4}}
\put(5,0){\line(0,1){1}} \put(6,0){\line(0,1){1}}
\put(7,0){\line(0,1){1}} \put(8,0){\line(0,1){1}}
\put(4.4,0.3){1} \put(5.4,0.3){2} \put(6.4,0.3){3}
\put(7.4,0.3){4}
\end{picture}

\begin{picture}(4,2.5)
\put(0,1){Region R2} \put(4,0){\line(0,1){2}}
\put(4,0){\line(1,0){2}} \put(4,1){\line(1,0){2}}
\put(4,2){\line(1,0){2}} \put(6,0){\line(0,1){2}}
\put(5,0){\line(0,1){2}} \put(4.4,1.3){1} \put(5.4,1.3){2}
\put(4.4,0.3){3} \put(5.4,0.3){4}
\end{picture}

\begin{picture}(4,2.5)
\put(0,1){Region R3}
\put(5,0){\line(1,0){2}}\put(4,1){\line(1,0){3}}
\put(4,1){\line(0,1){1}} \put(5,0){\line(0,1){2}}
\put(4,2){\line(1,0){2}} \put(6,0){\line(0,1){2}}
\put(7,0){\line(0,1){1}} \put(4.4,1.3){1} \put(5.4,1.3){2}
 \put(5.4,0.3){3} \put(6.4,0.3){4}
\end{picture}

\begin{picture}(4,2.5)
\put(0,1){Region R4} \put(4,0){\line(1,0){3}}
\put(4,1){\line(1,0){3}} \put(4,2){\line(1,0){1}}
\put(4,0){\line(0,1){2}} \put(5,0){\line(0,1){2}}
\put(6,0){\line(0,1){1}} \put(7,0){\line(0,1){1}}
\put(4.4,1.3){1}\put(4.4,0.3){2}
 \put(5.4,0.3){3} \put(6.4,0.3){4}
\end{picture}

\begin{picture}(4,2.5)
\put(0,1){Region R5} \put(4,0){\line(1,0){3}}
\put(4,1){\line(1,0){3}} \put(5,2){\line(1,0){1}}
\put(4,0){\line(0,1){1}} \put(5,0){\line(0,1){2}}
\put(6,0){\line(0,1){2}} \put(7,0){\line(0,1){1}}
\put(5.4,1.3){1}\put(4.4,0.3){2}
 \put(5.4,0.3){3} \put(6.4,0.3){4}
\end{picture}

Five customer arrival rates $\lambda$ are considered, with
$\lambda=3,6,9,12,15$ customers per unit time. There are eight customer types that capture a combination of location (1 to 4) and customer severity, high-level and low-level. Half of all customers at each location are high-level. All distances are measured from the center of one location to another. The distance between the centers of two neighboring locations is set to 1.0, and the distances between each pair of locations is computed as the Euclidean distance, resulting in distances between zero and three. The average service times capture travel time, which is assumed to be equal to the distance (zero to three), and service provided at the scene, which is assumed to be 12. Therefore, the average service times are between 12 and 15 for each customer and ambulance pair. Ambulance service rates are identical for high-level and low-level customers at the same location. 

\begin{figure}[h]
  \begin{center}
  \includegraphics{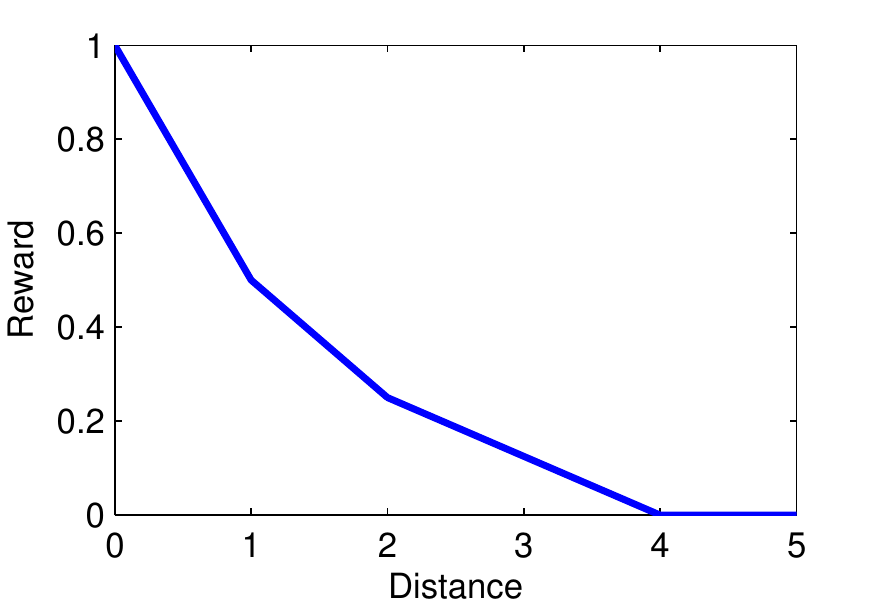}
\caption{Rewards for assigning an ambulance to a high-level customer
based on distance between the ambulance and the customer}
    \label{fig:reward}
  \end{center}
\end{figure}

The rewards are defined as a convex, piecewise linear function of
the distance between the ambulance and the customer. The rewards
for sending an ambulance to high-level customers that are 0, 1,
2, and 3 units away are 1, 1/2, 1/4, and 1/8, respectively. The
rewards for sending an ambulance to a low-level customer are 1/8
of the rewards for sending the same ambulance to a high-level
customer at the same location. Figure \ref{fig:reward}
illustrates the rewards for high-level customers. See \cite{yoon2018expected} for a discussion of how to weigh utilities of low-level customers relative to high-level customers.

There are five cases for the conditional arrival probabilities
(labeled C1 -- C5) to illustrate the effect of (geographic) customer
dispersion on the optimal solutions when a customer arrives. They are summarized in Table \ref{table:arrivalprob}. Case C1 corresponds to customers being uniformly dispersed throughout the region. Cases 2 and 3 correspond to customers being concentrated in one area either at the edge or in the center of the region, respectively. Cases 4 and 5 correspond to customers being concentrated in two of the areas at the edge or in the center of the region, respectively. With the five levels of $\lambda$, five regions, and five cases, there are 125 total scenarios.

\begin{table}\centering\small \caption{Conditional arrival probabilities of a
customer arriving at location $i$, $i=1,2,3,4$, given that a
customer has
arrived}\label{table:arrivalprob}\footnotesize\begin{tabular}{|cc|cccc|}
  \hline
Region & Arrival Case & $i=1$ & $i=2$ & $i=3$ & $i=4$ \\ \hline

\multirow{5}{*}{R1} & C1 & 0.25 & 0.25 & 0.25 & 0.25 \\
   & C2 & 0.7  & 0.1  & 0.1  & 0.1  \\
   & C3 & 0.1  & 0.7  & 0.1  & 0.1  \\
   & C4 & 0.4  & 0.1  & 0.1  & 0.4  \\
   & C5 & 0.1  & 0.4  & 0.4  & 0.1 \\ \hline
\multirow{5}{*}{R2} & C1 & 0.25 & 0.25 & 0.25 & 0.25 \\
   & C2 & 0.7  & 0.1  & 0.1  & 0.1  \\
   & C3 & 0.2  & 0.4  & 0.2  & 0.2  \\
   & C4 & 0.4  & 0.1  & 0.1  & 0.4  \\
   & C5 & 0.4  & 0.4  & 0.1  & 0.1 \\ \hline
\multirow{5}{*}{R3} & C1 & 0.25 & 0.25 & 0.25 & 0.25 \\
   & C2 & 0.7  & 0.1  & 0.1  & 0.1  \\
   & C3 & 0.1  & 0.7  & 0.1  & 0.1  \\
   & C4 & 0.4  & 0.1  & 0.1  & 0.4  \\
   & C5 & 0.1  & 0.4  & 0.4  & 0.1 \\ \hline
\multirow{5}{*}{R4} & C1 & 0.25 & 0.25 & 0.25 & 0.25 \\
   & C2 & 0.7  & 0.1  & 0.1  & 0.1  \\
   & C3 & 0.1  & 0.7  & 0.1  & 0.1  \\
   & C4 & 0.4  & 0.1  & 0.1  & 0.4  \\
   & C5 & 0.1  & 0.4  & 0.4  & 0.1 \\ \hline
\multirow{5}{*}{R5} & C1 & 0.25 & 0.25 & 0.25 & 0.25 \\
   & C2 & 0.7  & 0.1  & 0.1  & 0.1  \\
   & C3 & 0.1  & 0.1  & 0.7  & 0.1  \\
   & C4 & 0.1  & 0.4  & 0.1  & 0.4  \\
   & C5 & 0.4  & 0.1  & 0.4  & 0.1 \\ \hline
\end{tabular}\end{table}

All LP and MIP models are solved using Gurobi 4.5.1 on an Intel
Core 2 Duo 1.86 GHz processor with 2.94 GB of RAM. The
unrestricted model has 6,673 variables and 5,626 constraints.
Model PL has 6,801 variables and 9146 constraints. Gurobi
identified optimal solutions within 10 to 320 seconds of
execution time. We also compared the solutions to those of the heuristic that always sends the closest ambulance of those available. \added{Let $Z^H$ denote the solution value of the heuristic.}


First, the Model PL solutions are compared to the unrestricted model
solutions to identify the scenarios under which the optimal
unrestricted solution does not conform to a priority list. 
The optimal unrestricted solutions (to Model U) yield a priority list in 83 of the 125 scenarios, when priority lists are not explicitly enforced. In these 83 scenarios, corresponding Model PL and Model U solutions have the same objective function values. In the other 42 scenarios, the Model PL solution values are all within 0.66\% of the of the corresponding Model U solutions with an average gap of 0.12\%. Since the objective function values between corresponding Model PL and Model U solutions are extremely close, we focus on differences in the underlying solutions.
Table \ref{table:summaryPL} summarizes the number of
scenarios whose unrestricted policies are not priority list policies
(out of five total scenarios for each region and case). Regions R2 and R5 and cases C2, C3,
and C4 less frequently have optimal unrestricted policies that
conform to priority lists. 
In contrast, the heuristic of sending the closest ambulance often differs from the Model PL and U solutions. Table \ref{table:heuristic}(a) summarizes the average gaps between the closest ambulance heuristic solution values and the optimal Model U solution values across cases C1 -- C5. \added{Similarly, Table \ref{table:heuristic}(b) summarizes the average gaps between the closest ambulance heuristic solution values and the optimal Model PL solution values}. They report that the heuristic solution values are typically within 2.5\% of the optimal Model U and Model PL solution values, and that the differences increase with the customer arrival rate $\lambda$. 

\begin{table}\centering\small \caption{Number of unrestricted solutions to Model U that do not conform to priority lists}
\label{table:summaryPL}\footnotesize\begin{tabular}{|c|ccccc|c|}
  \hline
  & \multicolumn{5}{c|}{Case} & \\
Region & C1 & C2 & C3 & C4 & C5 & Total \\ \hline

R1 & 1 & 1 & 0 & 0 & 2 & 4 (of 25) \\
R2 & 0 & 5 & 0 & 5 & 1 & 11 (of 25) \\
R3 & 0 & 1 & 2 & 0 & 0 & 3 (of 25) \\
R4 & 1 & 3 & 1 & 1 & 0 & 6 (of 25) \\
R5 & 4 & 5 & 5 & 4 & 0 & 18 (of 25) \\ \hline Total & 6 (of 25)
& 15 (of 25) & 8 (of 25) & 10 (of 25) & 3 (of 25) & \textbf{42
(of 125)}
\\ \hline
\end{tabular}\end{table}



\begin{table}\centering\small \caption{\added{Average gaps (\%) between the closest ambulance heuristic solution values and the Model U and PL solution values}}
\label{table:heuristic}\footnotesize\begin{tabular}{|c|ccccc|c|}
 \multicolumn{7}{c}{(a) $(Z^U-Z^{PL})/Z^U$} \\
  \hline
Region & $\lambda=3$ & $\lambda=6$ & $\lambda=9$ & $\lambda=12$ & $\lambda=15$ & Average \\ \hline
R1 & 0.22\% & 1.36\% & 2.15\% & 2.44\% & 2.61\% & 1.76\% \\
R2 & 0.30\% & 1.10\% & 1.55\% & 1.75\% & 1.75\% & 1.29\% \\
R3 & 0.23\% & 1.35\% & 1.98\% & 2.21\% & 2.31\% & 1.62\% \\
R4 & 0.22\% & 1.33\% & 1.97\% & 2.35\% & 2.32\% & 1.64\% \\
R5 & 0.15\% & 1.26\% & 1.88\% & 2.16\% & 2.33\% & 1.56\% \\ \hline 
Average
& 0.23\% & 1.28\% & 1.90\% & 2.18\% & 2.27\% &\\ \hline
 \multicolumn{7}{c}{} \\
 \multicolumn{7}{c}{(b) $(Z^{PL}-Z^{H})/Z^{PL}$} \\
  \hline
  & \multicolumn{5}{c|}{$\lambda$} & \\
Region & $\lambda=3$ & $\lambda=6$ & $\lambda=9$ & $\lambda=12$ & $\lambda=15$ & Average \\ \hline
R1 & 0.41\% & 1.36\% & 2.15\% & 2.47\% & 2.61\% & 1.80\% \\
R2 & 0.44\% & 1.10\% & 1.57\% & 1.78\% & 1.79\% & 1.34\% \\
R3 & 0.26\% & 1.34\% & 1.98\% & 2.24\% & 2.31\% & 1.63\% \\
R4 & 0.37\% & 1.32\% & 2.01\% & 2.37\% & 2.40\% & 1.69\% \\
R5 & 0.26\% & 1.28\% & 1.92\% & 2.20\% & 2.37\% & 1.60\% \\ \hline 
Average
& 0.35\% & 1.28\% & 1.93\% & 2.21\% & 2.30\% &\\ \hline
\end{tabular}\end{table}

Next, we compare the solutions across different models focusing on
solutions to Model U, Model PL, and the closest ambulance heuristic.
Many of the scenarios have symmetric regions within them, and the
solutions to these scenarios are not unique in terms of the relative
distances between ambulances sent to customers. To address this issue,
policies are said to be the same if they dispatch symmetric ambulances.
For example, consider a customer arriving in $i=1$ in R2. Policies
that dispatch ambulance 2 or 3 are considered the same due to
symmetry.

\begin{figure}
    \centering
    \begin{subfigure}{0.5\linewidth}
    \centering 
     \fbox{\includegraphics[width = 0.8\textwidth]{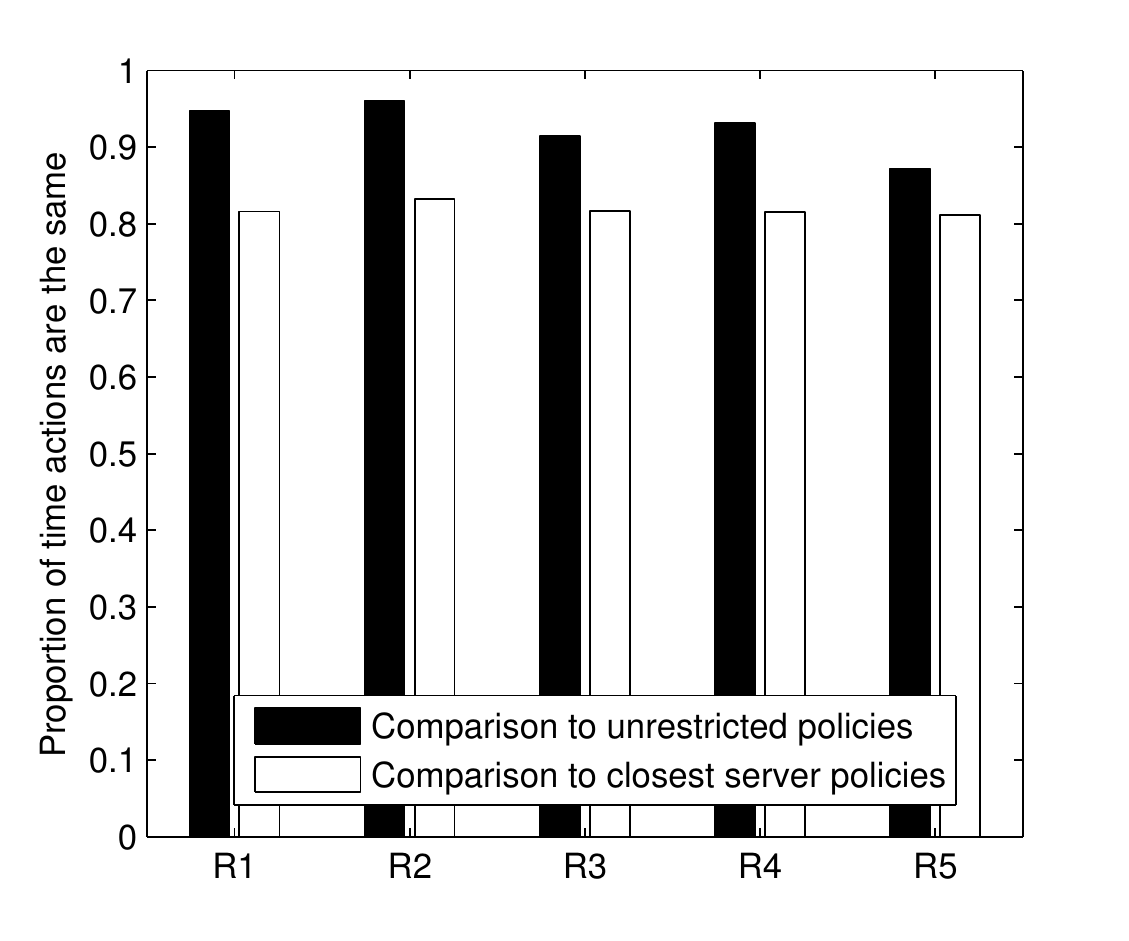}}
    \caption{by region}
    \end{subfigure} \hfill

    \begin{subfigure}{0.5\linewidth}
    \centering 
     \fbox{\includegraphics[width = 0.8\textwidth]{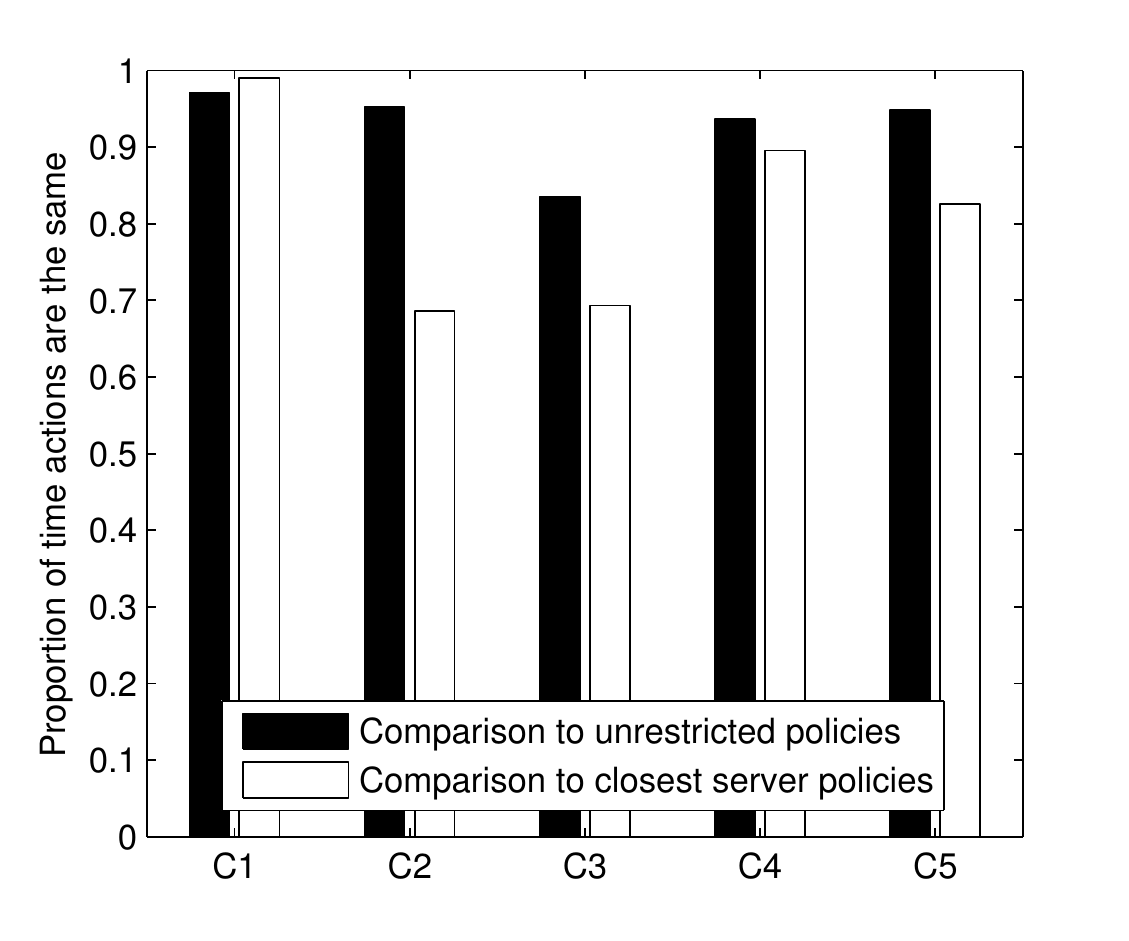}}
    \caption{by case}
    \end{subfigure} \hfill

    \begin{subfigure}{0.5\linewidth}
    \centering 
     \fbox{\includegraphics[width = 0.8\textwidth]{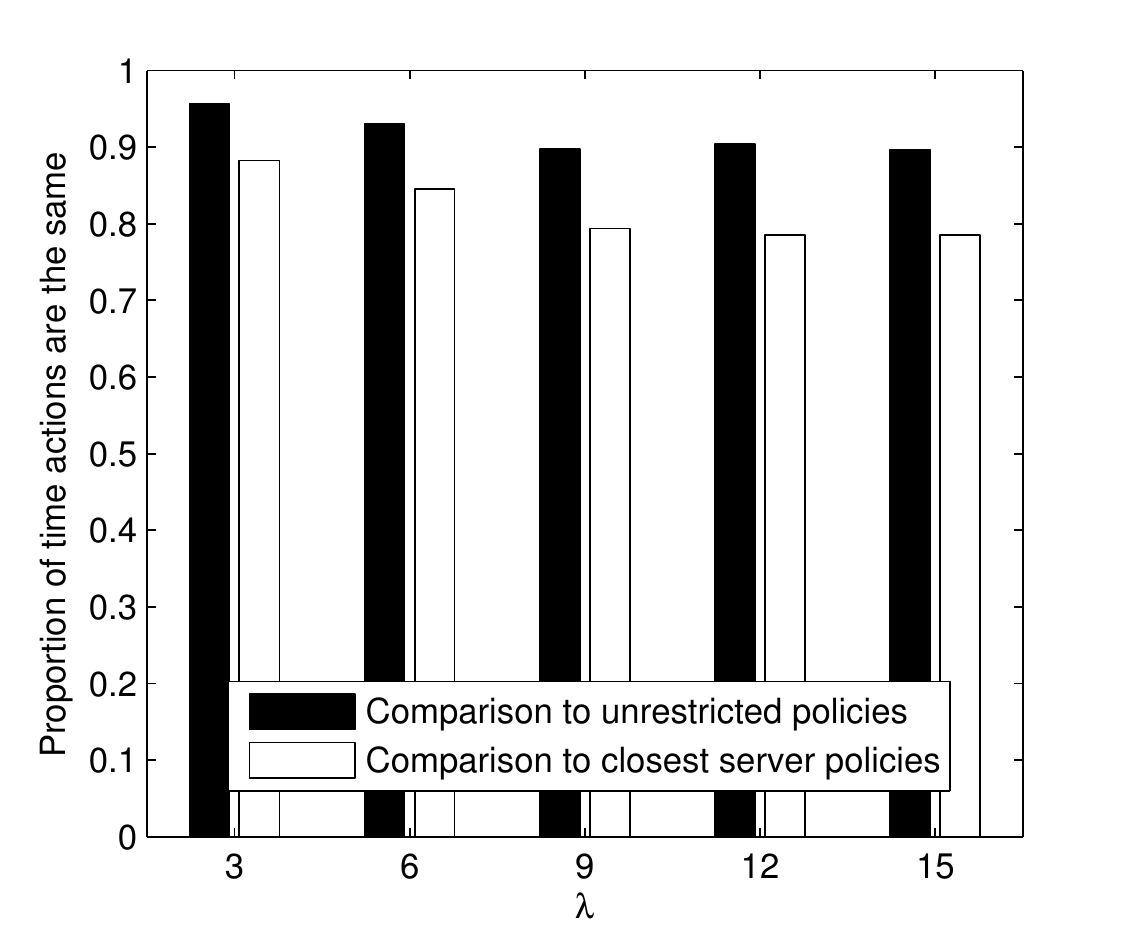}}
    \caption{by customer arrival rate}
    \end{subfigure}

    \caption{Comparisons of the Model PL solutions to the corresponding Model U and closest ambulance (server) heuristic solutions in terms of the proportion of time the same action is selected.}
    \label{fig:csr}
\end{figure}

Solutions to Model U often do not conform to priority lists, and when the Models U and PL solutions differ, the differences in the policies can be substantial. Figure \ref{fig:csr} provides three graphical comparisons between the optimal priority list and unrestricted solutions (the black bars) and between the optimal priority list and closest ambulance heuristic solutions (the white bars).
The comparisons between the priority list and unrestricted solutions shed light on whether Model PL results in any meaningful changes to the dispatching policies as compared to Model U,
and the comparisons between the priority list and closest ambulance heuristic solutions shed light on whether enforcing a priority list would lead to a simple policy that could more easily by obtained by a heuristic.

Each of the three subfigures in Figure \ref{fig:csr} portrays the proportion of time two of the policies are the same in the states where at least two actions are available to capture the degree of unintuitiveness.  Figure \ref{fig:csr}(a) suggests that policies are the same nearly the same proportion of time across each of the five regions, which suggests that the shape of the region may not greatly affect the policy. A different conclusion can be drawn when looking at Figure \ref{fig:csr}(b), which illustrates the results according to case. Here, the priority list solutions differ from the unrestricted and closest ambulance solutions the most for Cases 2 and 3. Finally, Figure \ref{fig:csr}(c) shows the similarities between the solutions as a function of $\lambda.$ The priority list solutions increasingly differ from the unrestricted and closest ambulance solutions as $\lambda$ increases, which is consistent with the observations in Table \ref{table:heuristic}. 
Both types of comparisons made in Figure \ref{fig:csr} suggest that Model PL identifies dispatching policies that would not be easily obtained without using the MIP modeling approach in many situations. 


Table  \ref{table:PL} reports the optimal priority lists for the
scenarios with region R5 and case C2 across different values of
$\lambda$. Each customer type is a combination of location (1 -- 4) and level (high (H) or low (L)). Recall that 70\% of the arriving customers are in
location 1. Each column reports the order that the four ambulances are sent to one of the eight types of customers. The first ambulance (see the rows with Priority $p=1$) indicates that the closest ambulance is always sent to each type of customer except to low-level customers at location 1.
Here, ambulance 1 is rationed for high-level calls. While the
differences between the unrestricted solutions are not reported here
for brevity, we note that the differences between the unrestricted
and priority list policies are in the backup coverage (i.e., the
priority 2 -- 4 ambulances to send to each type of customer). As
$\lambda$ increases, the need to ration ambulance 1 becomes apparent,
and when $\lambda=15$, ambulance 1 is the least preferred ambulance to
send to all customers except for high-level customers located at
1. The priority lists for high-level customer types match those of the heuristic of sending the closest ambulance across all levels of $\lambda$. The policies are markedly different for low-level customer types. 
Finally, note that the Model PL solutions are identical to
the closest ambulance heuristic solutions in only 3 of the 125 total scenarios,
which suggests that Model PL is valuable for improving upon heuristic policies as well as for providing an upper bound on the solution value for heuristics that conform to priority list policies.
In sum, the results from Table \ref{table:PL} and Figure \ref{fig:csr} demonstrate that Model PL leads to policies that are not easily obtained by using a simple heuristic.

\begin{table}\centering\small	\caption{Model PL priority lists for high-level (H) and low-level (L) customer types at each location for Region R5 and Case C2}																		
\label{table:PL}\footnotesize\begin{tabular}{|cc|cccccccc|}																			
\hline																			
&	& 	\multicolumn{8}{c|}{Customer type (Location, level)}		\\								$\lambda$	&	Priority ($p$)	&	1H	&	2H	&	3H	&	4H	
&	1L	& 	2L	& 	3L	& 	4L	 \\ \hline
\multirow{4}{*}{3}	&	1	&	1	&	2	&	3	&	4	&	3	&	2	&	3	&	4	 \\
	&	2	&	3	&	3	&	4	&	3	&	4	&	4	&	2	&	3	 \\
	&	3	&	4	&	1	&	1	&	1	&	1	&	3	&	1	&	2	 \\
	&	4	&	2	&	4	&	2	&	2	&	2	&	1	&	
4	&	1	 \\ \hline
\multirow{4}{*}{6}	&	1	&	1	&	2	&	3	&	4	&	3	&	2	&	3	&	4	 \\
	&	2	&	3	&	3	&	2	&	3	&	2	&	3	&	4	&	3	 \\
	& 	3	& 	2	& 	1	& 	4	& 	1	& 	4	&	4	&	2	&	2	 \\
	&	4	&	4	&	4	&	1	&	2	&	1	&	1	& 	
1	& 	1	 \\ \hline
\multirow{4}{*}{9}	&	1	&	1	&	2	&	3	&	4	&	2	&	2	&	3	&	4	 \\
	&	2	&	3	&	3	&	4	&	3	&	4	&	3	&	2	&	3	 \\
	&	3	&	4	&	4	&	2	&	2	&	3	&	4	&	4	&	2	 \\
	&	4	&	2	&	1	&	1	&	1	&	1	&	1	&	
1	&	1	 \\ \hline
\multirow{4}{*}{12}	&	1	&	1	&	2	&	3	&	4	&	2	&	2	&	3	&	4	 \\
	&	2	&	3	&	3	&	4	&	3	&	4	&	3	&	2	&	3	 \\
	&	3	&	2	&	4	&	2	&	2	&	3	&	4	&	4	&	2	 \\
	&	4	&	4	&	1	&	1	&	1	&	1	&	1	&	
1	&	1	 \\ \hline
\multirow{4}{*}{15}	&	1	&	1	&	2	&	3	&	4	&	2	&	2	&	3	&	4	 \\
	&	2	&	3	&	3	&	2	&	3	&	4	&	3	&	4	&	3	 \\
	&	3	&	2	&	4	&	4	&	2	&	3	&	4	&	2	&	2	 \\
	&	4	&	4	&	1	&	1	&	1	&	1	&	1	&	
1	&	1	 \\ \hline
\end{tabular}\end{table}	

\added{The models we introduce in this paper allows us to identify intuitive policies that are also optimal. One of the limitations of the models is that the state space grows exponentially with the number of servers, which limits the ability to solve large problem instances to optimality. To explore this issue of scalability, Table \ref{tab:cputime} reports the size of the state space, the number of variables, the number of constraints, and the mean run times (in seconds) for problems with different sizes. All instances were run on an Intel Core i7-1065G7 CPU at 1.3 GHz with 32 GB of RAM with a time limit of 7200 seconds. Table \ref{tab:cputime} indicates that instances with $m \ge 5$ typically require more than two hours to solve to optimality, which motivates the need for identifying explainable heuristics that closely mimic the optimal policies.
}

\begin{table}\centering\small \caption{\added{Mean run times for problems of different sizes with R5, C2, and $n=4$}}\label{tab:cputime}
\label{table:heuristic2}\footnotesize\begin{tabular}{|cc|ccc|ccc|}
  \hline
Number of & Size of & Number of  & Number of & Model PL & Number of  & Number of & Model U \\ 
servers $m$ & state space $|S|$ & constraints & variables & time (s) & constraints & variables & time (s) \\ \hline
3 & 125  & 1126   & 1237   & $<1$ & 1366   & 1309  & $<1$ \\
4 & 625  & 5626   & 6673   & 2.168& 9146   & 6801 & 55.6 \\
5 & 7776 & 28,126 & 36,317 & 86.9 & 67,246 & 36,517 & 6786* \\ \hline
	\multicolumn{8}{l}{* Time limit of 7200 s reached for 5 of 6 instances} \\
\end{tabular}\end{table}

\added{As noted earlier, we chose a larger state space with $|S|=(n+1)^m$ to highlight the value of the methodology introduced in this paper at the expense of scalability. In many settings, it is reasonable to reduce the state space to $|S|=2^m$ by assuming a constant service rate for each server (see \cite{yoon2021dynamic} for a discussion of this topic). With constant service times, there is a smaller state space as well as fewer variables and constraints. Table \ref{tab:cputime2} reports run times (in seconds) for problems with different numbers of servers $m$ under the constant service rate assumption. It indicates that solving problems with $m=10$ to $12$ is possible, which is realistic for many settings.
}

\begin{table}\centering\small \caption{\added{Run times for problems of different sizes with a constant service rate }}\label{tab:cputime2}
\label{table:heuristic2}\footnotesize\begin{tabular}{|cc|ccc|ccc|}
  \hline
Number of & Size of & Number of  & Number of & Model PL & Number of  & Number of & Model U \\ 
servers $m$ & state space $|S|$ & constraints & variables & time (s) & constraints & variables & time (s) \\ \hline
4  & 16   & 49     & 82     & $<1$ & 121    & 114   & $<1$ \\
5  & 32   & 95     & 194    & $<1$ & 375    & 244   & $<1$ \\
6  & 64   & 193    & 450    & $<1$ & 1093   & 522   & $<1$ \\
7  & 128  & 385    & 1026   & $<1$ & 2989   & 1124  & $<1$ \\
8  & 256  & 769    & 2306   & $<1$ & 7825   & 2434  & 2.1 \\
9  & 512  & 1537   & 5122   & 1.72 & 19,825 & 5284  & 1.2 \\
10 & 1024 & 3073   & 11,266 & 8.79 & 48,973 & 11,466& 8.9 \\
11 & 2048 & 6145   & 24,578 & 57.7 & 118,565& 24,820& 45.8 \\
12 & 4096 & 12,289  & 53,250 & 37.0 & 282,361& 53,538& 56.2 \\ 
13 & 8192 & 24,577   & 114,690    & 433 & 597,618   & 115,028 & 7200* \\ \hline
	\multicolumn{8}{l}{* Time limit of 7200 s reached} \\
\end{tabular}\end{table}

\subsection{Analysis of idling policies}
Next, we study the impact of allowing for idling policies.  In Model PLI, the number of state-action pairs increases with the addition of the action to idle. As a result, the corresponding Model U has 9,625 variables and 5,706 constraints, and Model PLI has has 9825 variables and 
20,970 constraints. 
Gurobi identified optimal solutions instances within 199 and 630 seconds of execution time for Models U and PLI, respectively.

We focus on identifying when idling is optimal, since many of the observations regarding the policies mirror those reported in the previous section. We report results for the five Region R5 and Case C2 scenarios.  When compared to
the corresponding scenario that does not allow idling, idling
improves the objective function values by 0.2\% when $\lambda=3$ to
17.7\% when $\lambda=15$. Table \ref{table:PLidle} reports the
optimal priority list policies when allowing idling for these five
scenarios. When $\lambda=3$, it is optimal to serve all customers
except in two cases when a single ambulance is available. An idling
policy occurs more frequently as $\lambda$ increases. Low-level
customers are served when all ambulances are available when $\lambda
\le 6$. It is optimal to serve some low-level customers across
all values of $\lambda$. As before, the need to ration ambulance 1 is apparent. In sum, this example illustrates how frequently low-level customers are ``lost'' due to idling when it is permissible as well as the differences in the policies as compared to when idling is not permitted (see Table \ref{table:PL}).
\begin{table}\centering\small\caption{Model PLI priority lists for high-level (H) and low-level (L) customer types at each location for for Region R5 and Case C2}
\label{table:PLidle}\small\begin{tabular}{|cc|cccccccc|}																		
\hline																	
&	& 	\multicolumn{8}{c|}{Customer type (Location, level)}	
\\																
$\lambda$	&	Priority ($p$)	&	1H	&	2H	&	3H	&	4H	
&	1L	& 	2L	& 	3L	& 	4L	 \\ \hline
\multirow{4}{*}{3}	&	1	&	1	&	2	&	3	&	4	&	1	&	2	&	3	&	4	 \\
	&	2	&	3	&	3	&	2	&	3	&	3	&	3	&	4	&	3	 \\
	&	3	&	4	&	4	&	1	&	1	&	4	&	4	&	2	&	2	 \\
&	4	&	2	&	1	&	4	&	2	&	Idle	&	1	
&	Idle	&	1	 \\ \hline
\multirow{4}{*}{6}	&	1	&	1	&	2	&	3	&	4	&	3	&	2	&	3	&	4	 \\
	&	2	&	3	&	3	&	2	&	3	&	Idle	&	Idle	&	4	&	Idle	 \\
	& 	3	& 	4	& 	4	& 	4	& 	2	& 	Idle	&	Idle	&	Idle	&	Idle	 \\
	&	4	&	2	&	1	&	1	&	1	&	Idle	&	
Idle	& 	Idle	& 	Idle	 \\ \hline
\multirow{4}{*}{9}	&	1	&	1	&	2	&	3	&	4	&	Idle	&	2	&	3	&	4	 \\
	&	2	&	3	&	3	&	4	&	3	&	Idle	&	Idle	&	Idle	&	Idle	 \\
	&	3	&	2	&	4	&	2	&	2	&	Idle	&	Idle	&	Idle	&	Idle	 \\
	&	4	&	4	&	1	&	1	&	1	&	Idle	&	
Idle	&	Idle	&	Idle	 \\ \hline
\multirow{4}{*}{12}	&	1	&	1	&	2	&	3	&	4	&	Idle	&	2	&	3	&	4	 \\
	&	2	&	3	&	3	&	4	&	3	&	Idle	&	Idle	&	Idle	&	Idle	 \\
	&	3	&	4	&	4	&	2	&	2	&	Idle	&	Idle	&	Idle	&	Idle	 \\
	&	4	&	2	&	1	&	1	&	1	&	Idle	&	
Idle	&	Idle	&	Idle	 \\ \hline
\multirow{4}{*}{15}	&	1	&	1	&	2	&	3	&	4	&	Idle	&	Idle	&	Idle	&	4	 \\
	&	2	&	3	&	3	&	4	&	3	&	Idle	&	Idle	&	Idle	&	Idle	 \\
	&	3	&	4	&	4	&	2	&	2	&	Idle	&	Idle	&	Idle	&	Idle	 \\
	&	4	&	2	&	1	&	1	&	1	&	Idle	&	
Idle	&	Idle	&	Idle	 \\ \hline
\end{tabular}\end{table}																			

\section{Conclusions}
Explainability is becoming increasingly important in the deployment of analytics-based decision support tools in the public sector \citep*{murdoch2019definitions} and is becoming a legal requirement in some countries \citep{goodman2017european}. Policies that are intuitive and transparent to public safety leaders are more likely to be trusted and used \citep*{levine2017new}.
As such, this paper is concerned with identifying intuitive ambulance dispatch policies using Markov decision process models that are optimal but not fully state-dependent.  To achieve this goal, we formulate new MIP models to identify optimal structured policies in a MDP that conform to priority lists.  The MIP models build upon the traditional linear programming formulation for solving MDP models. The MIP modeling framework introduces new binary variables for enforcing policies that conform to the desired structure without expanding the state space.
We introduce a base model that does not permit idling, and we extend this model to permit idling policies for low-level customers. In the models, a set
of constraints sets variable upper bounds for each of the MDP
actions. 

The MIP modeling approach is applied to a specific dispatching
problem where the priority list policies are desired. A series of
examples are examined to understand when the optimal
unrestricted solutions do not conform to priority lists; these
scenarios suggest the benefit of the proposed MIP models when the
calls are concentrated in a few areas in the region. A
significant empirical finding is that the MIP formulation does
not significantly degrade the solution values as compared to
the unrestricted solution values and improves upon the
solutions found with the heuristic of sending the closest
ambulance. It is worth noting that the MIP solutions result in
meaningful changes to the policies as compared to the
unrestricted policies, which suggest that it is not simple or
intuitive to alter the unrestricted policies to make them
conform to priority lists. Additional analysis and other metrics can justify whether simpler policies are actually beneficial.

A limitation of the approach in this paper is that the models do not scale well with the number of ambulances. Therefore, it is helpful to know when models that rely on simplifying assumptions and heuristics can be used. The empirical results suggest that in some settings it may be reasonable to make simplifying assumptions to identify priority list policies using approaches that are more tractable (e.g., \cite*{Ansari-etal-2012}). Note that if the service rates are the same for all customer types, an assumption that has often been shown to be approximately true in practice (see \cite*{Budge-etal-09}), the state space can be substantially reduced to solve larger problem instances.
Future work involves exploring explainable, data-driven dispatching models using, for example, approximate dynamic programming models that are more tractable yet offer public safety leaders and the public transparency. 

\section*{Acknowledgements}
This research was supported by the U.S. Department of the Army under
Grant Award Number W911NF-10-1-0176. The views and conclusions
contained in this document are those of the author and should not be
interpreted as necessarily representing the official policies,
either expressed or implied, of the U.S. Department of the Army.

\singlespacing
\bibliography{prioritylist}

\end{document}